%% file: NicKFU21.tex
\documentclass{ifacconf}

\usepackage{graphicx}       

\makeatletter
\def\endfigure{\end@float}
\def\endtable{\end@float}
\makeatother

\usepackage[caption=false]{subfig}
\usepackage{etoolbox}

\AtBeginEnvironment{figure}{\setlength{\captionwidth}{0.8\linewidth}}
\AtBeginEnvironment{table}{\setlength{\captionwidth}{0.8\linewidth}}
\graphicspath{{./img/}}
\usepackage{natbib}        

\usepackage{mathtools}
\usepackage{amssymb}
\usepackage{microtype}

\usepackage{xcolor}
\usepackage{xspace}

\usepackage{siunitx}
\usepackage{booktabs}

\usepackage{algorithm}
\usepackage[noend]{algpseudocode}

\algnewcommand\algorithmicforeach{\textbf{for each}}
\algdef{S}[FOR]{ForEach}[1]{\algorithmicforeach\ #1\ \algorithmicdo}

\usepackage{listings}
\definecolor{codegreen}{rgb}{0,0.6,0}
\definecolor{codegray}{rgb}{0.5,0.5,0.5}
\definecolor{codepurple}{rgb}{0.58,0,0.82}
\definecolor{backcolour}{rgb}{0.95,0.95,0.92}
\lstdefinestyle{mystyle}{
    backgroundcolor=\color{backcolour},
    commentstyle=\color{codegreen},
    keywordstyle=\color{magenta},
    numberstyle=\tiny\color{codegray},
    stringstyle=\color{codepurple},
    basicstyle=\ttfamily\scriptsize,
    breakatwhitespace=false,
    breaklines=true,
    captionpos=b,
    keepspaces=true,
    numbers=left,
    numbersep=5pt,
    showspaces=false,
    showstringspaces=false,
    showtabs=false,
    tabsize=2,
    xleftmargin = 2em,
}

\newcommand{\R}{\ensuremath\mathbb{R}}
\newcommand{\state}{\boldsymbol{x}}
\newcommand{\stateSpace}{\mathcal{X}}
\newcommand{\stateDim}{n}
\newcommand{\inpVar}{\boldsymbol{u}}
\newcommand{\optInpVar}{\boldsymbol{u}^{*}}
\newcommand{\inpVarSet}{\mathcal{U}}
\newcommand{\inpVarDim}{m}

\newcommand{\timeInt}{\mathbb{T}}
\newcommand{\flowVar}{\boldsymbol{\varphi}}

\newcommand{\weights}{\boldsymbol{\omega}}

\newcommand{\tddt}{\tfrac{\mathrm{d}}{\mathrm{d}t}}

\newcommand{\numData}{N_{\mathrm{data}}}
\newcommand{\numCP}{N_{\mathrm{phys}}}

\newcommand{\valNl}{4\xspace}
\newcommand{\valNn}{64\xspace}
\newcommand{\valNf}{20000\xspace}
\newcommand{\valNx}{100\xspace}
\newcommand{\valE}{800000\xspace}
\newcommand{\valMeanOCPst}{$3.65 \cdot 10^{-2}$~\si{\second}\xspace}
\newcommand{\valMAEalpha}{$7.53 \cdot 10^{-3}$~\si{\radian}\xspace}
\newcommand{\valMAEbeta}{$2.56 \cdot 10^{-2}$~\si{\radian}\xspace}

\newcommand{\bmat}[1]{\begin{bmatrix}#1\end{bmatrix}}
\newcommand{\bsmat}[1]{\begin{bsmallmatrix}#1\end{bsmallmatrix}}

\newcommand{\PINN}{PINN\xspace}
\newcommand{\DNN}{DNN\xspace}
\newcommand{\PINNs}{PINNs\xspace}
\newcommand{\DNNs}{DNNs\xspace}
\newcommand{\DAE}{DAE\xspace}
\newcommand{\MPC}{MPC\xspace}
\newcommand{\NMPC}{NMPC\xspace}

\newcommand{\ML}{ML\xspace}

\begin{document}


\begin{frontmatter}

\title{Physics-informed Neural Networks-based Model Predictive Control for Multi-link Manipulators\thanksref{footnoteinfo}}

\thanks[footnoteinfo]{The authors acknowledge funding from the DFG under Germany's Excellence Strategy -- EXC 2075 -- 390740016 and are thankful for support by the Stuttgart Center for Simulation Science (SimTech). }

\author[SC]{Jonas Nicodemus}
\author[ITM]{Jonas Kneifl}
\author[ITM]{J\"org Fehr}
\author[SC]{Benjamin Unger}

\address[SC]{Stuttgart Center for Simulation Science (SC SimTech), University of Stuttgart, Universit\"{a}tsstr.~32, 70569 Stuttgart, Germany (e-mail: \{jonas.nicodemus,benjamin.unger\}@simtech.uni-stuttgart.de)}
\address[ITM]{Institute of Engineering and Computational Mechanics, University of Stuttgart, Pfaffenwaldring 9, 70569 Stuttgart, Germany
 (e-mail: \{joerg.fehr,jonas.kneifl\}@itm.uni-stuttgart.de).}

\begin{abstract}                
We discuss \emph{nonlinear model predictive control} (\NMPC) for multi-body dynamics via physics-informed machine learning methods.
\emph{Physics-informed neural networks} (\PINNs) are a promising tool to approximate (partial) differential equations. 
\PINNs are not suited for control tasks in their original form since they are not designed to handle variable control actions or variable initial values.
We thus present the idea of enhancing \PINNs by adding control actions and initial conditions as additional network inputs. The high-dimensional input space is subsequently reduced via a sampling strategy and a zero-hold assumption. This strategy enables the controller design based on a \PINN as an approximation of the underlying system dynamics. The additional benefit is that the sensitivities are easily computed via automatic differentiation, thus leading to efficient gradient-based algorithms.
Finally, we present our results using our \PINN-based \MPC to solve a tracking problem for a complex mechanical system,  a multi-link manipulator.

\end{abstract}


\end{frontmatter}

\section{Introduction}

\emph{Model predictive control} (\MPC) is a flexible and intuitive control scheme that lets us impose constraints and helps to operate complex systems optimally. The major challenge for \MPC is the repetitive solution of an optimal control problem. Even with today's computing power, efficient model representation to be real-time capable remains the bottleneck. This issue is especially prominent in systems with a fast dynamic, like robotic manipulators, where operation speed relates to increased productivity. 

In this work, we study a tracking problem for a multi-link manipulator, see section~\ref{sec:problem_formulation} for further details, with an a-priori unknown tracking trajectory. In this scenario, standard linearization strategies for the nonlinear dynamics, which are used to speed up the computation, cannot be implemented without further challenges. Instead, a repetitive numerical evaluation of the underlying nonlinear dynamics is required. Since the time required by standard time integration schemes may pose a critical constraint during the solution of the optimal control problem, we propose replacing the time-integration with a \emph{machine learning}~(\ML) approach.  Since standard \ML techniques typically require an extensive training data set and cannot compete with state-of-the-art time-integrators \citep{OtnGBPPSZ21}, we propose a physics-informed approach, see \cite{KarKLPWY21} for a recent overview, thus exploiting the underlying known physical law during the training process.  More precisely, we replace the nonlinear dynamics with a \emph{physics-informed neural network} (\PINN), initially introduced by \cite{RaiPK19}. In addition, we  exploit the efficient computation of partial derivatives of the \emph{deep neural network} (\DNN) via automatic differentiation, see for instance \cite{BayPRS18}.

Our main results are the following:
\begin{enumerate}
	\item Following ideas presented by~\cite{AntCSSJH21}, we detail in section~\ref{subsec:pinn_based_mpc} how to replace the nonlinear dynamics with a \PINN approximation in the context of  \emph{nonlinear model predictive control} (\NMPC). Due to the efficient computation of the partial derivative of the \PINN with respect to the controller, the optimal control problem can be solved efficiently with a gradient-based method, without any adjoint computations (for gradient-based methods) or acceleration strategies (for nonlinear programming methods).
	\item In section~\ref{sec:results} we apply the strategy to a PowerCube serial robot (cf.~\cite{FehSSE20}) and demonstrate that  replacing the numerical time-integration of the nonlinear dynamics with a \PINN speeds up the computation time while retaining a sufficient accuracy within the \NMPC framework.
\end{enumerate}

The use of ML techniques in \MPC is not new. For instance, in \cite{AkeT06} and \cite{HerKTA18} the controller itself is replaced by a neural network. In the context of a-priori unknown tracking trajectories, such an approach is either not feasible or requires learning in a high-dimensional space. Instead, similarly as in this work, the nonlinear dynamics are replaced by a recurrent neural network in \cite{WuRGC21} or a \PINN in \cite{ArnK21,AntCSSJH21}. 
However, \cite{ArnK21} use a different strategy for dealing with variable control inputs, and compared to \cite{AntCSSJH21}, we test our \PINN-based \MPC on a mechanical system that features faster dynamics than their examples.

Let us briefly outline the structure of the paper.
First, we postulate the precise problem setting, including the PowerCube serial robot as sample application in section~\ref{sec:problem_formulation}.  A short review of \NMPC and \PINNs in sections~\ref{subsec:NMPC} and~\ref{subsec:PINN}, respectively, is followed by the \PINN-based \NMPC framework in section~\ref{subsec:pinn_based_mpc}. We demonstrate the efficiency of the framework with a numerical example of multi-link manipulator in section~\ref{sec:results}.

\paragraph*{Notation}
For a positive semidefinite matrix $\boldsymbol{Q}\in \R^{n \times n}$, we define the weighted (euclidean) seminorm
\begin{displaymath}
	\|\cdot\|_{\boldsymbol{Q}}\colon \R^n\to \R_{\geq 0},\qquad \state \mapsto \sqrt{\state^\top \boldsymbol{Q}\state}.
\end{displaymath}

\section{Problem formulation}
\label{sec:problem_formulation}
On the time-interval $\timeInt \subseteq \R$ we study control systems of the form
\begin{subequations}
	\label{eqn:controlSystem}
	\begin{align}
		\label{eqn:controlSystem:equation}
		\dot{\state}(t) &= \boldsymbol{f}(\state(t),\inpVar(t)),\\
		\label{eqn:controlSystem:initialCondition}
		\state(t_0) &= \state_0,
    \end{align}
\end{subequations}
for some initial time $t_0\in\timeInt$ and $\state\colon \timeInt \to \stateSpace \subseteq \R^{\stateDim}$,  $\state_0\in\R^{\stateDim}$, $\inpVar\colon \timeInt \to \inpVarSet\subseteq\R^{\inpVarDim}$ the \emph{state}, \emph{initial value} and \emph{control}, respectively.
We assume that $f$ is continuous and (locally) Lipschitz-continuous with respect to the state, such that the initial value problem~\eqref{eqn:controlSystem} has a unique (weak) solution for each $\inpVar\in L^\infty(\timeInt,\inpVarSet)$, cf.~\citet[Thm.~54]{Son98}. The corresponding solution operator (also called \emph{flow}) that maps the control and initial value to the solution at time $t\geq t_0$, is denoted with
\begin{equation}
	\label{eqn:flowMap}
	\flowVar(t,\inpVar,\state_0) = \state(t).
\end{equation}
For our further presentation it will be beneficial to view the control system~\eqref{eqn:controlSystem:equation} as the under-determined \emph{differential-algebraic equation} (\DAE)
\begin{equation}
	\label{eqn:DAEformulation}
	0 = \boldsymbol{F}(\state(t),\dot{\state}(t),\inpVar(t))\vcentcolon=  \dot{\state}(t) - \boldsymbol{f}(\state(t),\inpVar(t)).
\end{equation}
Let us emphasize that, in principle, we could start with a descriptor system instead of the control system~\eqref{eqn:controlSystem}. However, for the sake of simplicity, we restrict ourselves to control systems of the form~\eqref{eqn:controlSystem}.

\begin{figure}
	\centering
 	\def\svgwidth{.95\linewidth}
    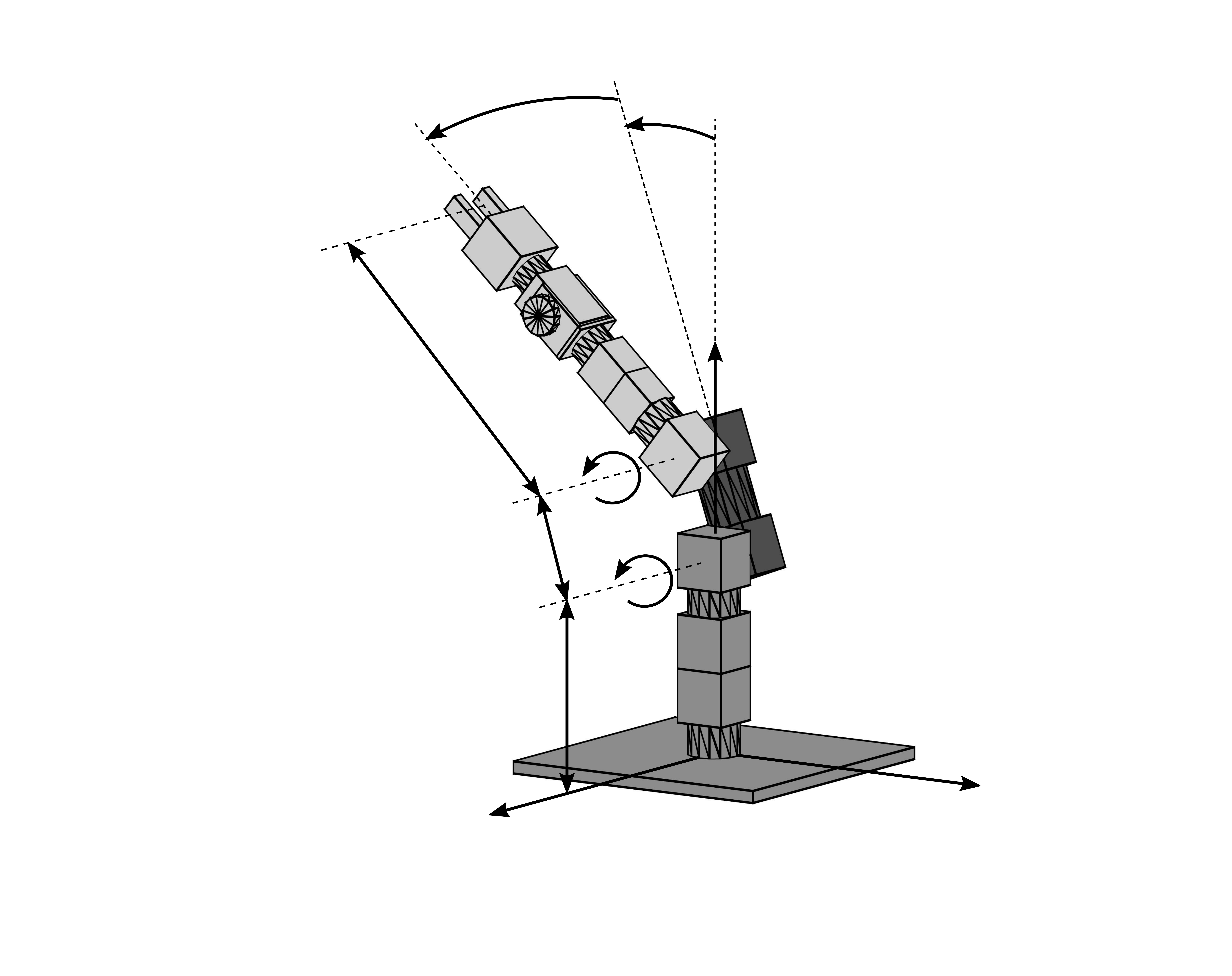
	\caption{Sketch of the PowerCube serial robot.}
	\label{fig:robot}
\end{figure}

As an exemplary real-world application, we study a multi-link manipulator consisting of Schunk PowerCubes modules, as described in \cite{FehSSE20}.
Multi-link manipulators play an important role in the automation process of the manufacturing industry and are a heavily studied research topic in fields like control theory, see for instance \cite{SpoHV20}.
The control problem of a multi-link manipulator is simple but still challenging.
Thanks to previous works \citep{Kar20}, the dynamical model of the manipulator is already identified and simulation data is available for comparison.
The schematic sketch of the here studied manipulator is shown in Fig.~\ref{fig:robot}.

The resulting equations of motion can be read as 
\begin{equation*}
    \label{eqn:motion}
    \boldsymbol{M}(\boldsymbol{q}) \ddot{\boldsymbol{q}} + \boldsymbol{k}(\boldsymbol{q}, \dot{\boldsymbol{q}}) = \boldsymbol{h}(\boldsymbol{q}, \dot{\boldsymbol{q}}) + \boldsymbol{B} \inpVar,
\end{equation*}
with the generalized coordinates $\boldsymbol{q}=[\alpha, \beta]^\top \in \R^2$, the nonsingular mass matrix $\boldsymbol{M}(\boldsymbol{q}) \in \R^{2 \times 2}$, the vector of centrifugal, coriolis and gyroscopic forces $\boldsymbol{k}(\boldsymbol{q}, \dot{\boldsymbol{q}}) \in \R^{2}$, the vector of applied forces $\boldsymbol{h}(\boldsymbol{q}, \dot{\boldsymbol{q}}) \in \R^{2}$, the input matrix $\boldsymbol{B} \in \R^{2 \times 2}$ and the input vector $\inpVar \in \R^2$, which consists of the motor currents of the A and B revolute joint.
Introducing $\state\vcentcolon=[\boldsymbol{q}, \dot{\boldsymbol{q}}]^\top$ yields the first-order reformulation
\begin{equation*}
    \boldsymbol{0} = \bmat{\boldsymbol{I} & \boldsymbol{0} \\ \boldsymbol{0} & \boldsymbol{M}(\boldsymbol{q})}\bmat{\dot{\boldsymbol{q}}\\\ddot{\boldsymbol{q}}} - 
    \bmat{\dot{\boldsymbol{q}} \\ \boldsymbol{h}(\boldsymbol{q}, \dot{\boldsymbol{q}}) - \boldsymbol{k}(\boldsymbol{q}, \dot{\boldsymbol{q}})} -
    \bmat{0 \\ \boldsymbol{B} \inpVar},
\end{equation*}
which is in the form of \eqref{eqn:DAEformulation}, with $\state \in \R^{4}$.

\section{Methods}
\label{sec:methods}
The practical implementation of a controller on a digital control unit typically requires a temporal discretization of the continuous-time control system~\eqref{eqn:controlSystem}.  We follow \cite{GruP11} and introduce the time-grid
\begin{displaymath}
	t_0 < t_1 < \ldots < t_N.
\end{displaymath}
For simplicity we choose an equidistant time grid, i.e., $t_k = k\tau + t_0$ with constant sampling period $\tau>0$.  Introducing the shifted input $\boldsymbol{u}_k \vcentcolon= \boldsymbol{u}(\cdot+k\tau)$ for $k=0,1,\ldots,N-1$ we thus obtain
\begin{equation}
	\label{eqn:controlSystem:discrete}
	\state_{k+1} = \flowVar(\tau,\inpVar_k,\state_k)
\end{equation}
with $\state_k = \state(t_k)$ for $k=0,1,\ldots,N$. If we make an additional zero-hold assumption for the control input, i.e., $\tddt \inpVar_{k|(0,\tau)} \equiv 0$, then~\eqref{eqn:controlSystem:discrete} resembles a discrete-time control system. In this case, by abuse of the notation, we use the symbol $\inpVar_k$ both to denote the control function and its constant value on the interval $(0,\tau)$.

\subsection{Nonlinear Model Predictive Control}
\label{subsec:NMPC}

Suppose we have a control system in the form of~\eqref{eqn:controlSystem:discrete}, with the state measured at discrete time instants $t_k$ as described above, then the tracking problem is to find suitable control inputs $\inpVar_k$ such that $\state_k$ follows a given reference $\state^{\mathrm{ref}}_k$ as good as possible. 
This problem can be solved by applying \NMPC.

For this purpose, we introduce a cost function $\ell\colon \stateSpace \times \stateSpace \times \inpVarSet \to \R_{\ge 0}$. In the quadratic case,  this function may be chosen as
\begin{equation}
    \label{eqn:stage_costs}
    \ell(\state^{\mathrm{ref}}_k, \state_k, \inpVar_k) = \| \state^{\mathrm{ref}}_k - \state_k \|^{2}_{Q} + \| \inpVar_k \|^{2}_{R}.
\end{equation}

\begin{figure}
	\centering
    \includegraphics[width=0.95\linewidth]{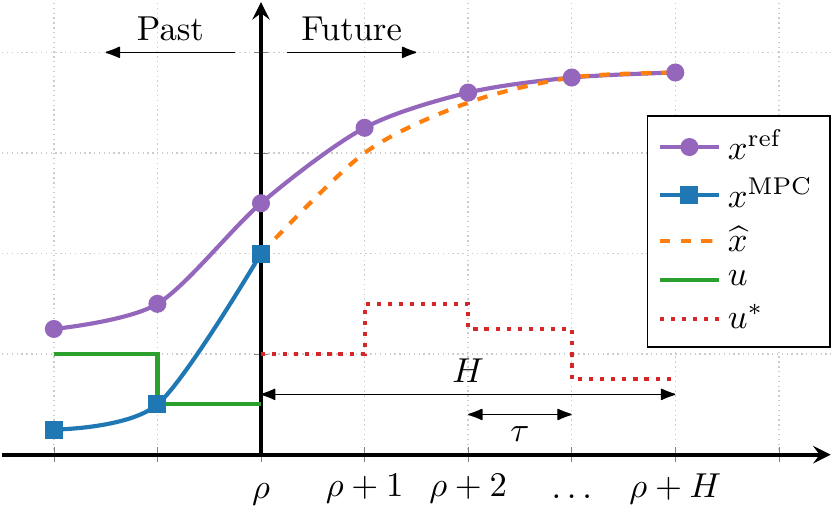}
	\caption{Basic \MPC scheme for the horizon width $H=4$.}
	\label{fig:MPC_scheme}
\end{figure}

Following the \NMPC scheme, visualized in Fig.~\ref{fig:MPC_scheme}, at each time instant $k = \rho$ the discrete-time optimal control problem for the moving horizon with horizon width $H$ for the current time instant $\rho$ is given by
\begin{subequations}
\label{eqn:disc_ocp}
\begin{align}
    \min_{\inpVar_{\rho},\ldots, {\inpVar_{\rho + H - 1}}} & \sum_{k=\rho}^{\rho + H - 1} \ell(\state^{\mathrm{ref}}_k, \state_k, \inpVar_k) \label{eqn:disc_ocp:costs}\\
    \text{s.\,t.} \qquad & \state_{k+1} = \flowVar(\tau,\inpVar_k,\state_k),\\
    & \inpVar_k \in \inpVarSet, \state_k \in \stateSpace.
\end{align}
\end{subequations}
Then, the optimal control input sequence $\optInpVar_k$ for $k= \rho, \ldots, \rho + H - 1$ for the current time instant $\rho$ is obtained, by solving \eqref{eqn:disc_ocp}.
From this sequence the first control input~$\optInpVar_\rho$ is applied to the system until $t_{\rho + 1}$, then the procedure repeats.
The optimal control problem \eqref{eqn:disc_ocp} is solved numerically either: (i) directly, via the Hamiltonian-Jacobi-Bellmann equation, (ii) via gradient-based methods, where the gradients are computed via the adjoint method, or as a (iii) nonlinear programming problem with methods such as sequential quadratic programming or the interior-point method, see for instance~\cite{NocW06}. 

In many applications, $\flowVar$ is not available and a numerical approximation~$\widehat{\flowVar}$ is used, for instance obtained via numerical time-integration. In the case of the tracking problem with a-priori known tracking trajectory, one can linearize the nonlinear system in an offline-step to speed up the computation in the online stage. This strategy is used for instance in~\cite{FehSSE20} for the optimal control of the present robot manipulator. We emphasize that the reference trajectory needs to be known a-priori and the linearization introduces an additional source of error. We thus cannot apply this strategy and restrict ourselves to the fully nonlinear case.

\subsection{Physics-Informed Neural Networks}
\label{subsec:PINN}
\begin{figure}
	\centering
	\includegraphics[width=.95\linewidth]{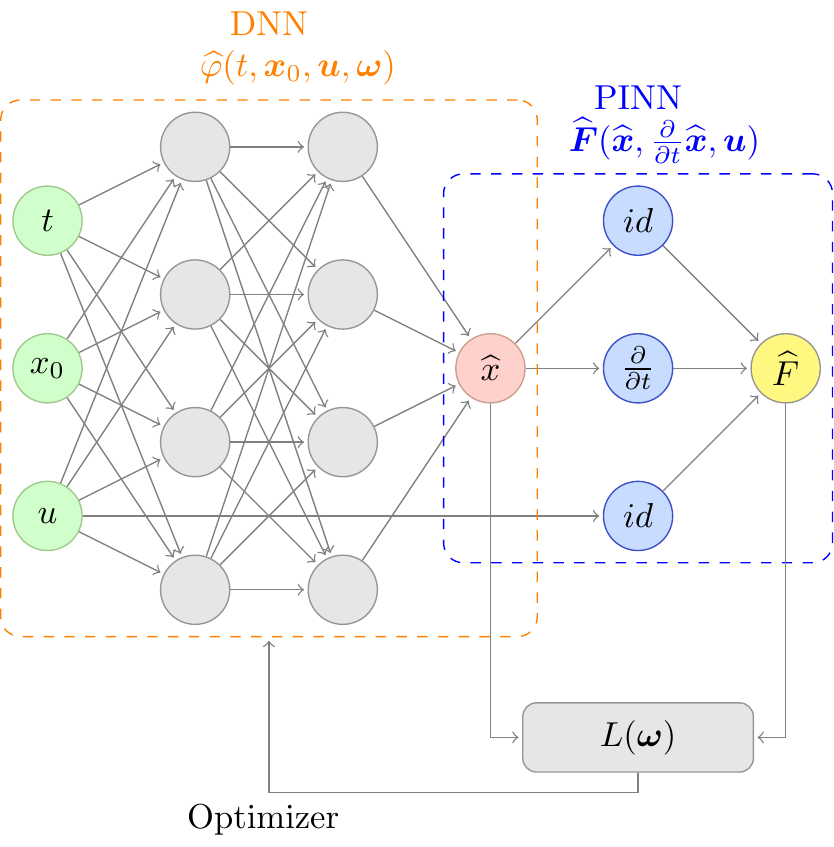}
	\caption{\PINN architecture.}
	\label{fig:PiNN}
\end{figure}

If various solution trajectories of~\eqref{eqn:controlSystem} are available, then \DNNs can be used to approximate the solution map~\eqref{eqn:flowMap} via a supervised learning task. In more detail, a (data-based) loss-function based on the available solution trajectories is minimized with respect to so-called \emph{weights} $\weights\in\R^{p}$ such that the \DNN $\widehat{\flowVar}$ approximates $\flowVar$, i.e.,
\begin{displaymath}
	\widehat{\state}(t) \vcentcolon= \widehat{\flowVar}(t,\state_0,\inpVar,\weights) \approx \flowVar(t,\state_0,\inpVar),
\end{displaymath}
for all admissible $(t,\state_0,\inpVar)$. One strategy to achieve this, is to use existing data
\begin{displaymath}
	\state_i \vcentcolon= \flowVar(t_i,\state_{0,i},\inpVar_i), \qquad i=1,\ldots,\numData
\end{displaymath}
and minimize the mean-squared error loss-function
\begin{equation*}
    L_{\mathrm{data}}(\weights) \vcentcolon= \frac{1}{\numData} \sum_{i=1}^{\numData}
    \|\widehat{\flowVar}(t_i,\state_{0,i},\inpVar_i,\weights) - \state_i\|^2.
\end{equation*}
The idea of \PINNs, introduced in \cite{RaiPK19}, is to add the differential equation to the loss function, as shown in Fig.~\ref{fig:PiNN}, to robustify the network further and allow for a good approximation even in a data-poor regime.
The actual \PINN results from inserting the approximated flow $\widehat{\flowVar}$ and its derivative $\tfrac{\partial}{\partial t} \widehat{\flowVar}$, which can be computed via automatic differentiation, into the differential equation~\eqref{eqn:DAEformulation}, yielding the residual
\begin{equation*}
    \widehat{\boldsymbol{F}}(t, \state_0, \inpVar, \weights) \vcentcolon= \boldsymbol{F}(\widehat{\flowVar}(t,\state_{0},\inpVar,\weights), \tfrac{\partial}{\partial t} \widehat{\flowVar}(t,\state_{0},\inpVar,\weights), \inpVar).
\end{equation*}
Note that this approach requires sufficiently smooth activation functions in the \DNN, for instance the hyperbolic tangent.
The loss function corresponding to the differential equation is then given by
\begin{equation*}
    L_{\mathrm{phys}}(\weights) \vcentcolon= \frac{1}{\numCP} \sum_{i=1}^{\numCP}
    \|\widehat{\boldsymbol{F}}(t_i, \state_{0,i}, \inpVar_i, \weights)\|^2,
\end{equation*}
where the residual is evaluated on a finite set of so-called \emph{collocation points} $(t_i,\state_{0,i},\inpVar_i)$.
Let us emphasize that the collocation points can be chosen arbitrarily without the need to solve the control system~\eqref{eqn:controlSystem}.

To ensure that both the data and the differential equation are approximated sufficiently well, we minimize the \DNN with respect to the combined loss function.
\begin{equation*}
    L(\weights)= L_{\mathrm{data}}(\weights) + L_{\mathrm{phys}}(\weights).
\end{equation*}
In practical applications, it may be important to introduce a weighting of the two components of the loss function. For notational convenience, we assume that such a scaling is already implicitly encoded in~$\widehat{\boldsymbol{F}}$.

\begin{figure*}
    \centering
    \subfloat[\PINN in self-loop prediction.]{
        \centering
        \includegraphics{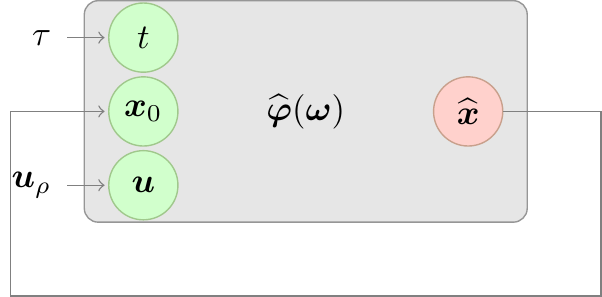}
        \label{fig:pinn_self_loop}
    }
    \hfill
    \subfloat[\PINN-based \MPC in closed-loop.]{
        \centering
        \includegraphics{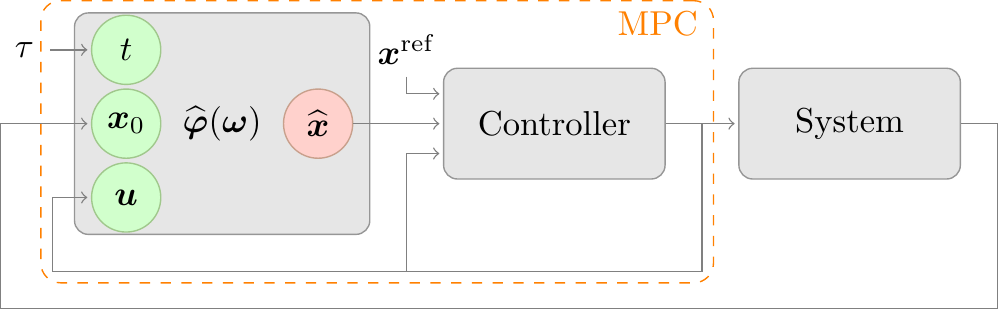}
        \label{fig:PINN_based_MPC}
    }
    \caption{On the left \PINN in self-loop prediction and on the right \PINN-based \MPC connected to a system.}
    \label{fig:PINN_loops}
\end{figure*}

\subsection{\PINN-based \MPC} 
\label{subsec:pinn_based_mpc}
In principle, we could implement the \PINN approach as presented in the previous section. Nevertheless, the specific control application results in two major challenges: First, to obtain accurate approximations, we have to sample the infinite-dimensional input space $L^\infty(\timeInt,\inpVarSet)$. Second, while the \PINN approximation is inherently smooth, the solution of the control problem~\eqref{eqn:controlSystem} is typically not differentiable with respect to time.  Although the universal approximation theorem guarantees arbitrarily accurate approximations of continuous functions, this normally comes with the cost of a large number of neurons.

To remedy these issues within the \NMPC framework, we make the following observation: If $\boldsymbol{f}$ is sufficiently smooth and the zero-hold assumption is applied, then also $\state_{|(k\tau,(k+1)\tau)}$ is smooth. Moreover, within the interval $(k\tau,(k+1)\tau)$, the zero-hold assumption reduces the infinite-dimensional input space to the finite-dimension set $\inpVarSet\subseteq \R^{\inpVarDim}$. We immediately arrive at the following strategy. Given an initial value $\state_k$ and values $\inpVar_k$ for the control, train a \PINN that is able to predict $\state_{k+1}$. Note that in order to exploit the underlying dynamics in the training process, we additionally need an explicit dependency on the time variable. We conclude that the zero-hold assumption allows us to reduce the sampling space from $\timeInt\times \stateSpace\times L^\infty(\timeInt,\inpVarSet)$ to $[0,\tau]\times \stateSpace \times \inpVarSet$, thus rendering the learning task feasible. 

Two remarks are in order: First, the initial value is included in the sampling space, thus deviating from the original \PINN approach presented by \cite{RaiPK19}, where the networks are trained for a fixed initial value. We emphasize that this is not due to our particular approach but inherent to the \MPC framework. Second, our network is not restricted to predictions over the interval $[0,\tau]$.  To compute predictions for $t>\tau$, we simply take the \PINN approximation $\widehat{\flowVar}(\tau,\state_0,\inpVar_0)$ at time $t=\tau$ as new initial value for the time interval $[\tau,2\tau]$ and repeat this process iteratively. This approach, which we refer to as self-loop prediction, is illustrated in Fig.~\ref{fig:pinn_self_loop}. To corresponding \MPC scheme where the nonlinear dynamics are replaced with the \PINN approximation is presented in Fig.~\ref{fig:PINN_based_MPC}.

\section{Results}
\label{sec:results}
We test our methodology introduced in section~\ref{sec:methods} with the multi-link manipulator robot model discussed in section~\ref{sec:problem_formulation}. Our approach is implemented in TensorFlow\footnote{https://www.tensorflow.org/}. The network topology, including all relevant hyperparameters, is discussed in section~\ref{subsec:numericalResults:networkTopology}. To ensure that the \PINN approximation is sufficiently accurate, we first present in section~\ref{subsec:numericalResults:selfLoop} the prediction for a given input sequence. There the \PINN operates in self-loop prediction mode as described in section~\ref{subsec:pinn_based_mpc}. The tracking problem for the multi-link manipulator is then solved via the \PINN-based \MPC in closed-loop (cf.~section~\ref{subsec:numericalResults:closedLoop}). In the following, we use the symbols $\widehat{\alpha}$, $\widehat{\beta}$, $\alpha_{\mathrm{MPC}}$,  and $\beta_{\mathrm{MPC}}$ to denote the \PINN predictions for the angles $\alpha$ and $\beta$ from Fig.~\ref{fig:robot},  and the resulting angles from the \MPC scheme, respectively.

To ensure reproducibility of the conducted experiments, the code for the numerical examples is publicly available under the \texttt{doi:10.5281/zenodo.5520662}. 
All simulations are performed on an Apple M1 chip.

\subsection{Network topology and training}
\label{subsec:numericalResults:networkTopology}

Our network consists of \valNl layers, each hidden layer is activated by the hyperbolic tangent and contains \valNn neurons.
For the data-based loss function $L_{\mathrm{data}}$, we use $\numData = \valNx$ data points, which only contain a subset of possible initial values, thus not a single solution of the nonlinear system~\eqref{eqn:controlSystem} is required, see also \cite{RaiPK19}. 

Our numerical experiments have shown that it is reasonable to train the network for a slightly larger sampling period than it will be evaluated later on, i.e., we train the network for the time interval $[0,\tilde{\tau}]$ with $\widetilde{\tau} \vcentcolon= 0.25~\si{\second}>0.2~\si{\second} =\vcentcolon \tau$, and use
\begin{align*}
    \stateSpace &= [-\pi, \pi]^2 \times [-2.5, 2.5]^2, &
    \inpVarSet &= [-0.5, 0.5]^2.
\end{align*}
For the physics-informed loss function $L_{\mathrm{phys}}$, we sample $[0,\widetilde{\tau}] \times \stateSpace \times \inpVarSet \subseteq \R^{7}$ with $M_{\mathrm{phys}} = \valNf$ collocation points, generated via Latin Hypercube sampling \citep{McKBC79}.
All hyper parameter are chosen as a result of a grid search.
Let us emphasize that despite the small state dimension ($\stateDim=4$), we have to sample a 7-dimensional space, rendering this a challenging problem from an approximation perspective.

The network training, i.e., the minimization of the loss function, is performed via L-BFGS, a quasi-Newton, full-batch gradient-based optimization algorithm \citep{LiuN89} with $\valE$ iterations (referred to as epochs in the literature).
To optimize the the training results, we choose to generate new data and collocation after the first 400000 epochs.

\subsection{\PINN in self-loop prediction}
\label{subsec:numericalResults:selfLoop}

\begin{figure}
    \centering
    \includegraphics[trim=0 .4cm 0 .35cm, clip]{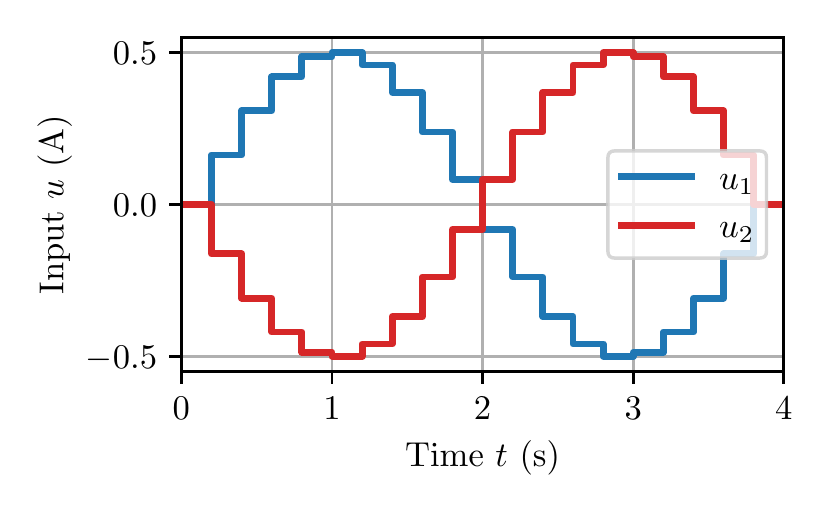}
    \caption{Testing input sequence $u$ for the open-loop application.}
    \label{fig:testing_u}
\end{figure}

\begin{figure}
    \centering
    \includegraphics[trim=0 .4cm 0 .35cm, clip]{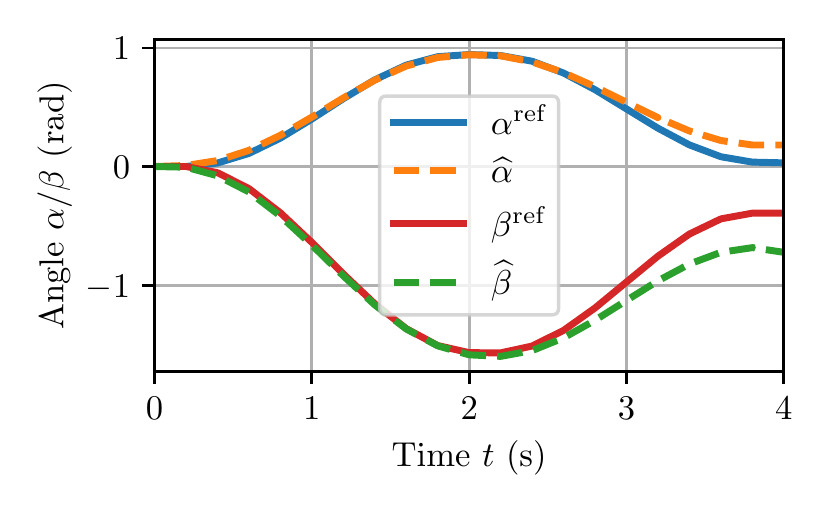}
    \caption{Open-loop simulation result from \PINN operated in self-loop prediction mode.}
    \label{fig:pinn_open_loop}
\end{figure}

To ensure a sufficiently good prediction capability of the \PINN approximation, we first run the \PINN in self-loop mode (cf.~Fig.~\ref{fig:pinn_self_loop}) with the testing input sequence presented in Fig.~\ref{fig:testing_u}. 
The sampling period is set to $\tau=0.2$~\si{\second}, thus this experiment results in 20 loop iterations.
The corresponding open-loop prediction is presented in Fig.~\ref{fig:pinn_open_loop}.
We observe that the discrepancy between the \PINN prediction and the nominal dynamics is relatively small for $t\leq 2$~\si{\second}. For $t\geq 2$~\si{\second}, the accumulation of the approximation errors due to the self-loop mode is noticeable and increases further over time. In the \MPC context, the solution of the optimal control problem is robust regarding to this error accumulation for larger time horizons (cf.~\cite{Sch21}), such that we expect the \PINN approximation quality for the dynamics is sufficient within the \MPC framework.

\subsection{\PINN in closed-loop}
\label{subsec:numericalResults:closedLoop}
\begin{figure}
    \centering
    \includegraphics[trim=0 .4cm 0 .7cm, clip, width=0.85\linewidth]{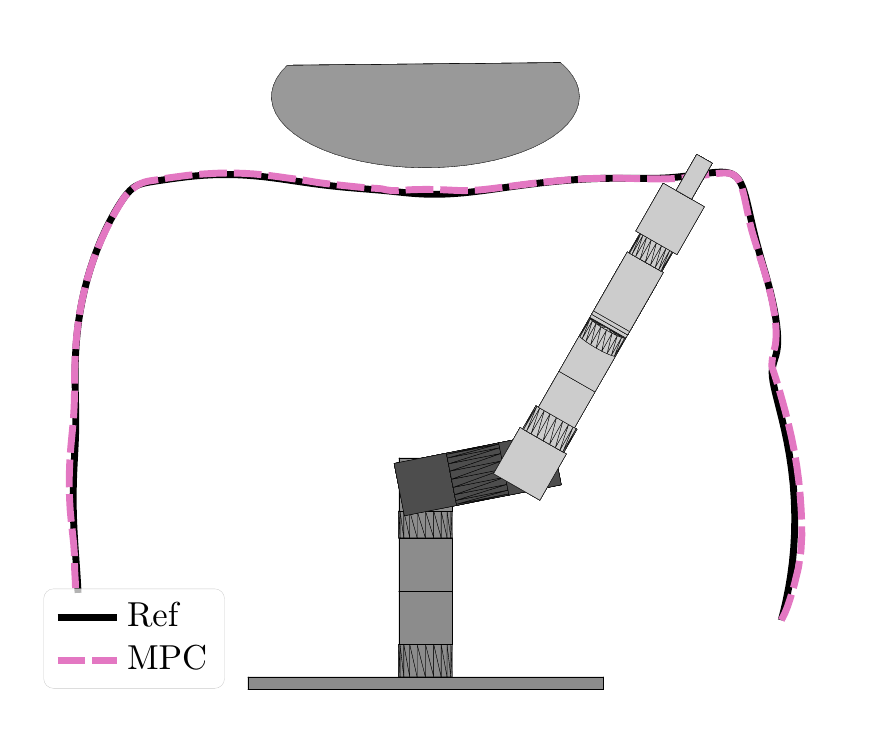}
    \caption{Reference trajectory and resulting \PINN-based \MPC trajectory for the tracking problem of the mulit-link manipulator displayed at time $t=8$\,\si{\second}.}
    \label{fig:testing_trajectory}
\end{figure}

We now apply the \PINN-based \MPC discussed in section~\ref{subsec:pinn_based_mpc} to a tracking problem for the multi-link manipulator (cf.~section~\ref{sec:problem_formulation}).
For this purpose,  we use the testing trajectory shown in Fig.~\ref{fig:testing_trajectory} and its corresponding reference states $\state^{\mathrm{ref}}$, taken from \cite{FehSSE20}.
The trajectory is motivated by an obstacle avoidance and includes motion reversals in both joints to investigate the control performance on friction phenomena caused by the harmonic drive gears that underly highly nonlinear friction effects.

The plant system is simulated by the Runge–Kutta–Fehlberg method (RK45), the \MPC sampling period is set to $\tau=0.2$\,\si{\second} and a horizon width of $H=5$ is used. 
These parameters are taken from \cite{FehSSE20}.
The \MPC cost function~\eqref{eqn:stage_costs} is parameterized by the diagonal matrices $\boldsymbol{Q}$ and $\boldsymbol{R}$.
We chose the matrices as
\begin{equation*}
    \boldsymbol{Q} \vcentcolon= \bsmat{1 & 0 & 0 & 0 \\ 
    0 & 1 & 0 & 0\\
    0 & 0 & 0 & 0\\
    0 & 0 & 0 & 0} , \qquad
    \boldsymbol{R} \vcentcolon= \bsmat{1 \cdot 10^{-6} & 0\\ 0 & 1 \cdot 10^{-6}},
\end{equation*}
which means we only control the positions of the angels.

\begin{figure}[hbt]
    \centering
    \includegraphics[trim=0 .4cm 0 .35cm, clip]{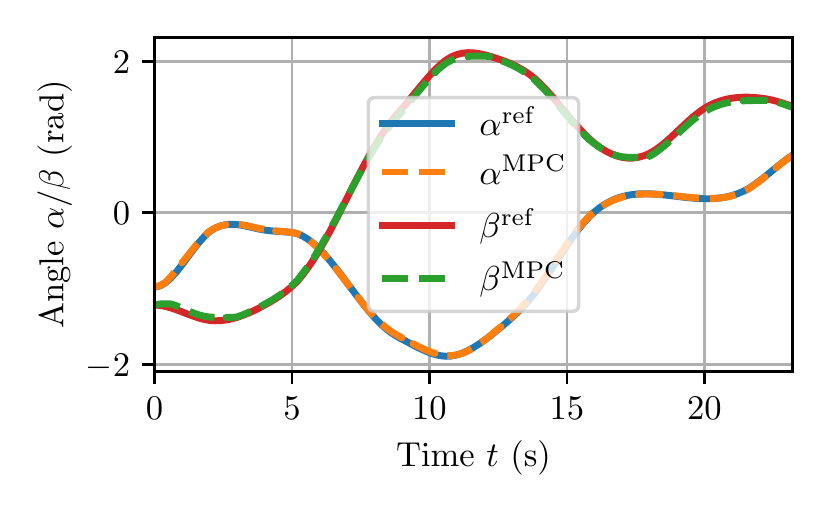}
    \caption{Solution of the \PINN-based \MPC for the tracking problem of the multi-link manipulator.}
    \label{fig:pinn_mpc}
\end{figure}
\begin{figure}[hbt]
    \centering
    \includegraphics[trim=0 .4cm 0 .35cm, clip]{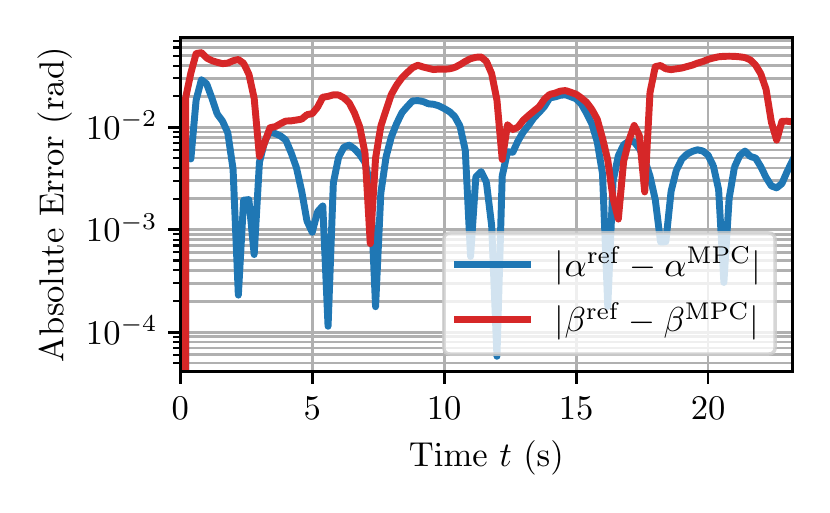}
    \caption{Absolute error over time of the solution of \PINN-based \MPC for the tracking problem of the multi-link manipulator.}
    \label{fig:mpc_error}
\end{figure}
\begin{figure}
    \centering
    \includegraphics[trim=0 .4cm 0 .35cm, clip]{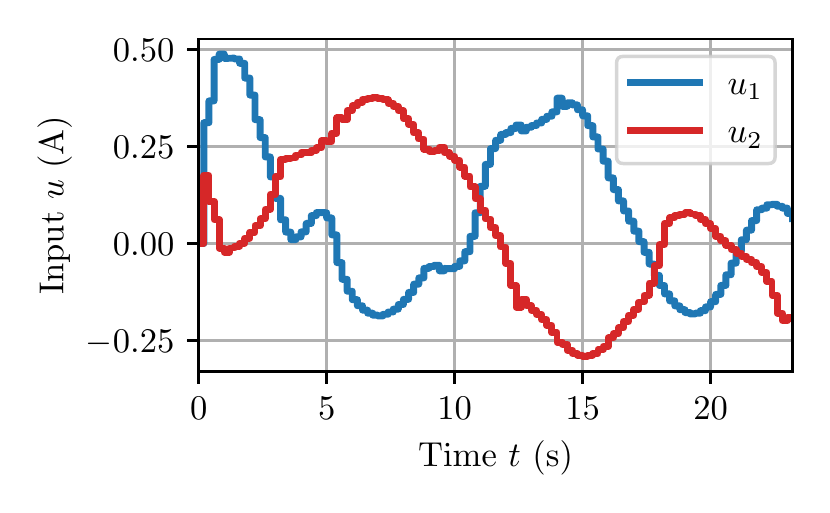}
    \caption{\PINN-based \MPC optimal $u$.}
    \label{fig:pinn_mpc_u}
\end{figure}
Figs.~\ref{fig:pinn_mpc} and~\ref{fig:mpc_error} illustrate that the \NMPC scheme with \PINN approximation for the nonlinear dynamics follows the reference trajectory closely with a mean absolute error of \valMAEalpha and \valMAEbeta for $\alpha$ and $\beta$, respectively. The computed control input is presented in Fig.~\ref{fig:pinn_mpc_u}. The solution of the \NMPC problem is computed in average in \valMeanOCPst, which is clearly below the sampling period.

Let us emphasize that already the computation of $\state_{k+1}$ from $\state_{k}$ is more expensive with a classical time-integration scheme such as RK45 or the explicit Euler method, than with the \PINN approximation, see Table~\ref{tb:execution_times}.
\begin{table}[htb]
\centering
\caption{Execution times for different Runge-Kutta methods compared with the \PINN, with $\tau=0.2$~\si{\second} and $h=0.02$~\si{\second} for Euler and RK4. }\label{tb:execution_times}
\begin{tabular}{lrrrr} 
	\toprule
    & \PINN & Euler & RK4 & RK45 \\\midrule
Mean (\si{\second}) & 4.14e-04 & 5.93e-04 & 2.35e-03 & 8.62e-03 \\
Median (\si{\second}) & 4.12e-04 & 5.90e-04 & 2.34e-03 & 1.77e-03 \\ \bottomrule
\end{tabular}
\end{table}

The main reason for the numerical integrators' performance is that these methods cannot use the sampling period as time-step, but need a smaller step-size $h$.
For instance, the explicit Euler method is unstable for this system with step-size $h = \tau$. Instead, we use $h=0.02$\,\si{\second} for the explicit Euler method, which in turn results in a slower execution time. Two additional remarks are in order: First, the evaluation of the \PINN is implemented in TensorFlow. 
If the network with final weights is directly implemented using state-of-the-art BLAS libraries, we expect an additional speedup. 
Second, the computation of the gradient of the nonlinear dynamics concerning the control can be performed in parallel to the evaluation of the \PINN (in contrast to adjoint-based methods), giving an additional speedup within the \NMPC framework.

\section{Conclusion}
We studied the applicability of a physics-informed machine learning approach, namely physics-informed neural networks (\PINNs) \citep{RaiPK19}, in the context of nonlinear model predictive control (\NMPC) for a multi-link robot manipulator.  Following ideas presented by \cite{AntCSSJH21}, we discussed how \PINNs could be used to replace the nonlinear dynamics with an efficient-to-evaluate surrogate model. Due to automatic differentiation, the \PINN approximation offers, in addition, a cheap and easy-to-implement computation of the partial derivative of the state with respect to the control input, thus paving the way for gradient-descent methods to solve the \NMPC problem without the need of further acceleration strategies.
In our numerical example for a tracking problem for a multi-link manipulator, we have shown that the \PINN-based approximation of the nonlinear dynamics outperforms classical time-integrators in terms of computational time while at the same time being accurate enough to solve the tracking problem.

\bibliography{literature}

\end{document}

%% file: img/schunk_scematic.pdf_tex
\begingroup%
  \makeatletter%
  \providecommand\color[2][]{%
    \errmessage{(Inkscape) Color is used for the text in Inkscape, but the package 'color.sty' is not loaded}%
    \renewcommand\color[2][]{}%
  }%
  \providecommand\transparent[1]{%
    \errmessage{(Inkscape) Transparency is used (non-zero) for the text in Inkscape, but the package 'transparent.sty' is not loaded}%
    \renewcommand\transparent[1]{}%
  }%
  \providecommand\rotatebox[2]{#2}%
  \newcommand*\fsize{\dimexpr\f@size pt\relax}%
  \newcommand*\lineheight[1]{\fontsize{\fsize}{#1\fsize}\selectfont}%
  \ifx\svgwidth\undefined%
    \setlength{\unitlength}{4202.78920091bp}%
    \ifx\svgscale\undefined%
      \relax%
    \else%
      \setlength{\unitlength}{\unitlength * \real{\svgscale}}%
    \fi%
  \else%
    \setlength{\unitlength}{\svgwidth}%
  \fi%
  \global\let\svgwidth\undefined%
  \global\let\svgscale\undefined%
  \makeatother%
  \begin{picture}(1,0.60226666)%
    \lineheight{1}%
    \setlength\tabcolsep{0pt}%
    \put(0,0){\includegraphics[width=\unitlength,trim=0 18cm 0 10cm, clip,page=1]{schunk_scematic.pdf}}%
    \put(0.71084423,0.35594414){\color[rgb]{0,0,0}\makebox(0,0)[lt]{\begin{minipage}{0.25464413\unitlength}\raggedright  \end{minipage}}}%
    \put(0.48044242,0.17026683){\color[rgb]{0,0,0}\makebox(0,0)[lt]{\begin{minipage}{0.12095595\unitlength}\raggedright $T_1$\end{minipage}}}%
    \put(0.43497803,0.34634669){\color[rgb]{0,0,0}\makebox(0,0)[lt]{\begin{minipage}{0.12095595\unitlength}\raggedright $T_2$\end{minipage}}}%
    \put(0.35603442,0.03335816){\color[rgb]{0,0,0}\makebox(0,0)[lt]{\begin{minipage}{0.14896683\unitlength}\raggedright $e_3$\end{minipage}}}%
    \put(0.79859515,0.0600958){\color[rgb]{0,0,0}\makebox(0,0)[lt]{\begin{minipage}{0.14896683\unitlength}\raggedright $e_1$\end{minipage}}}%
    \put(0.60032809,0.40369625){\color[rgb]{0,0,0}\makebox(0,0)[lt]{\begin{minipage}{0.14896683\unitlength}\raggedright $e_2$\end{minipage}}}%
    \put(0.42351511,0.13266937){\color[rgb]{0,0,0}\makebox(0,0)[lt]{\begin{minipage}{0.14896683\unitlength}\raggedright $\ell_0$\end{minipage}}}%
    \put(0.41078285,0.23505152){\color[rgb]{0,0,0}\makebox(0,0)[lt]{\begin{minipage}{0.14896683\unitlength}\raggedright $\ell_1$\end{minipage}}}%
    \put(0.33311641,0.36312232){\color[rgb]{0,0,0}\makebox(0,0)[lt]{\begin{minipage}{0.14896683\unitlength}\raggedright $\ell_2$\end{minipage}}}%
    \put(0.4358286,0.57090687){\color[rgb]{0,0,0}\makebox(0,0)[lt]{\begin{minipage}{0.14896683\unitlength}\raggedright $\beta$\end{minipage}}}%
    \put(0.53824875,0.54240952){\color[rgb]{0,0,0}\makebox(0,0)[lt]{\begin{minipage}{0.17825095\unitlength}\raggedright $\alpha$\end{minipage}}}%
    \put(0,0){\includegraphics[width=\unitlength, trim=0 18cm 0 10cm, clip,page=2]{schunk_scematic.pdf}}%
    \put(0.74257596,0.26144838){\color[rgb]{0,0,0}\makebox(0,0)[lt]{\begin{minipage}{0.22729316\unitlength}\raggedright $m_1, J_1$\end{minipage}}}%
    \put(0.09307083,0.41540019){\color[rgb]{0,0,0}\makebox(0,0)[lt]{\begin{minipage}{0.22729316\unitlength}\raggedright $m_2, J_2$\end{minipage}}}%
    \put(0,0){\includegraphics[width=\unitlength,trim=0 18cm 0 10cm, clip,page=3]{schunk_scematic.pdf}}%
    \put(0.81640081,0.44442081){\color[rgb]{0,0,0}\makebox(0,0)[lt]{\begin{minipage}{0.22729316\unitlength}\raggedright $g$\end{minipage}}}%
    \put(0.27671207,0.27263786){\color[rgb]{0,0,0}\makebox(0,0)[lt]{\begin{minipage}{0.13006048\unitlength}\raggedright joint B\end{minipage}}}%
    \put(0.2907476,0.18589495){\color[rgb]{0,0,0}\makebox(0,0)[lt]{\begin{minipage}{0.13006048\unitlength}\raggedright joint A\end{minipage}}}%
  \end{picture}%
\endgroup%

%% file: NicKFU21.bbl
\begin{thebibliography}{18}
\providecommand{\natexlab}[1]{#1}
\providecommand{\url}[1]{\texttt{#1}}
\providecommand{\urlprefix}{URL }
\expandafter\ifx\csname urlstyle\endcsname\relax
  \providecommand{\doi}[1]{doi:\discretionary{}{}{}#1}\else
  \providecommand{\doi}{doi:\discretionary{}{}{}\begingroup
  \urlstyle{rm}\Url}\fi

\bibitem[{Antonelo et~al.(2021)Antonelo, Camponogara, Seman, de~Souza,
  Jordanou, and Hubner}]{AntCSSJH21}
Antonelo, E.A., Camponogara, E., Seman, L.O., de~Souza, E.R., Jordanou, J.P.,
  and Hubner, J.F. (2021).
\newblock Physics-informed neural nets-based control.
\newblock \emph{arXiv preprint arXiv:2104.02556}.

\bibitem[{Arnold and King(2021)}]{ArnK21}
Arnold, F. and King, R. (2021).
\newblock State–space modeling for control based on physics-informed neural
  networks.
\newblock \emph{Engineering Applications of Artificial Intelligence}, 101,
  104195.
\newblock \doi{10.1016/j.engappai.2021.104195}.

\bibitem[{Baydin et~al.(2018)Baydin, Pearlmutter, Radul, and
  Siskind}]{BayPRS18}
Baydin, A.G., Pearlmutter, B.A., Radul, A.A., and Siskind, J.M. (2018).
\newblock Automatic differentiation in machine learning: a survey.
\newblock \emph{J. Mach. Learn. Res.}, 18.

\bibitem[{Fehr et~al.(2020)Fehr, Schmid, Schneider, and Eberhard}]{FehSSE20}
Fehr, J., Schmid, P., Schneider, G., and Eberhard, P. (2020).
\newblock Modeling, simulation, and vision-/{MPC}-based control of a
  {PowerCube} serial robot.
\newblock \emph{Applied Sciences}, 10(20), 7270.
\newblock \doi{10.3390/app10207270}.

\bibitem[{Gr{\"u}ne and Pannek(2011)}]{GruP11}
Gr{\"u}ne, L. and Pannek, J. (2011).
\newblock \emph{Nonlinear model predictive control}.
\newblock Springer, London.
\newblock \doi{10.1007/978-3-319-46024-6}.

\bibitem[{Hertneck et~al.(2018)Hertneck, K\"ohler, Trimpe, and
  Allg\"ower}]{HerKTA18}
Hertneck, M., K\"ohler, J., Trimpe, S., and Allg\"ower, F. (2018).
\newblock Learning an approximate model predictive controller with guarantees.
\newblock \emph{IEEE Control Systems Letters}, 2(3).
\newblock \doi{10.1109/LCSYS.2018.2843682}.

\bibitem[{Kargl(2020)}]{Kar20}
Kargl, A. (2020).
\newblock {DMD Methoden zur Identifikation von Reibkraftmodellen für einen
  Schunk PowerCube-Roboter}.
\newblock Bachelor thesis {BSC-120}, Institute of Engineering and Computational
  Mechanics, University Stuttgart.

\bibitem[{Karniadakis et~al.(2021)Karniadakis, Kevrekidis, Lu, Perdikaris,
  Wang, and Yang}]{KarKLPWY21}
Karniadakis, G.E., Kevrekidis, I.G., Lu, L., Perdikaris, P., Wang, S., and
  Yang, L. (2021).
\newblock Physics-informed machine learning.
\newblock \emph{Nature Reviews Physics}, 3, 422--440.
\newblock \doi{10.1038/s42254-021-00314-5}.

\bibitem[{Liu and Nocedal(1989)}]{LiuN89}
Liu, D.C. and Nocedal, J. (1989).
\newblock On the limited memory {BFGS} method for large scale optimization.
\newblock \emph{Math. Program.}, 45(1), 503--528.
\newblock \doi{10.1007/BF01589116}.

\bibitem[{McKay et~al.(1979)McKay, Beckman, and Conover}]{McKBC79}
McKay, M.D., Beckman, R.J., and Conover, W.J. (1979).
\newblock A comparison of three methods for selecting values of input variables
  in the analysis of output from a computer code.
\newblock \emph{Technometrics}, 42(1), 55--61.

\bibitem[{Nocedal and Wright(2006)}]{NocW06}
Nocedal, J. and Wright, S. (2006).
\newblock \emph{Numerical Optimization}.
\newblock Springer Series in Operations Research and Financial Engineering.
  Springer New York.

\bibitem[{Otness et~al.(2021)Otness, Gjoka, Bruna, Panozzo, Peherstorfer,
  Schneider, and Zorin}]{OtnGBPPSZ21}
Otness, K., Gjoka, A., Bruna, J., Panozzo, D., Peherstorfer, B., Schneider, T.,
  and Zorin, D. (2021).
\newblock An extensible benchmark suite for learning to simulate physical
  systems.
\newblock \emph{ArXiv e-print 2108.07799}.

\bibitem[{Raissi et~al.(2019)Raissi, Perdikaris, and Karniadakis}]{RaiPK19}
Raissi, M., Perdikaris, P., and Karniadakis, G. (2019).
\newblock Physics-informed neural networks: A deep learning framework for
  solving forward and inverse problems involving nonlinear partial differential
  equations.
\newblock \emph{J. Comput. Phys.}, 378, 686--707.
\newblock \doi{10.1016/j.jcp.2018.10.045}.

\bibitem[{\r{A}kesson and Toivonen(2006)}]{AkeT06}
\r{A}kesson, B.M. and Toivonen, H, T. (2006).
\newblock A neural network model predictive controller.
\newblock \emph{Journal of Process Control}, 16(9), 937--946.
\newblock \doi{10.1016/j.jprocont.2006.06.001}.

\bibitem[{Schaller(2021)}]{Sch21}
Schaller, M. (2021).
\newblock \emph{Sensitivity Analysis and Goal Oriented Error Estimation for
  Model Predictive Control}.
\newblock Ph.D. thesis, University of Bayreuth.

\bibitem[{Sonntag(1989)}]{Son98}
Sonntag, E.D. (1989).
\newblock \emph{Mathematical Control Theory}.
\newblock Springer, New York, 2 edition.
\newblock \doi{10.1007/978-1-4612-0577-7}.

\bibitem[{Spong et~al.(2020)Spong, Hutchinson, and Vidyasagar}]{SpoHV20}
Spong, M.W., Hutchinson, S., and Vidyasagar, M. (2020).
\newblock \emph{Robot modeling and control}.
\newblock John Wiley \& Sons.

\bibitem[{Wu et~al.(2021)Wu, Rincon, Gu, and Christofides}]{WuRGC21}
Wu, Z., Rincon, D., Gu, Q., and Christofides, P.D. (2021).
\newblock Statistical machine learning in model predictive control of nonlinear
  processes.
\newblock \emph{Mathematics}, 9(16), 1912.
\newblock \doi{10.3390/math9161912}.

\end{thebibliography}
